\documentclass{article}
\usepackage{graphicx}
\usepackage{amsmath,amssymb}
\newcommand{\ba}{\begin{array}}
\newcommand{\ea}{\end{array}}
\def\be{\begin{equation}}
\def\ee{\end{equation}}

\def\b{\beta}

\def\l{\lambda}

\def\f{\frac}

\def\q{\quad}
\def\be{\begin{equation}}
\def\ee{\end{equation}}
\def\bel{\begin{equation}\label}

\def\bb{\textbf{\textit{b}}}

\def\v{\textbf{\textit{v}}}

\def\y{\textbf{\textit{y}}}
\def\fb{\textbf{\textit{f}}}

\def\0{\textbf{\textit{0}}}
\def\L{{\mathcal L}}

\def\bmat{\left[\begin{matrix}}
\def\emat{\end{matrix}\right]}
\def\bpmat{\left(\begin{bmatrix}}
\def\epmat{\end{bmatrix}\right)}

\newtheorem{theorem}{Theorem}[section]

\newtheorem{remark}[theorem]{Remark}
\newtheorem{example}[theorem]{Example}

\begin{document}

\centerline{\Large\bf On the Method of Partial Fractions with Matrix }
\centerline{\Large\bf Coefficients and Applications}
\vspace{0.5cm}

\centerline{\large Ruben Airapetyan\footnote{email: rhayrape@kettering.edu}} 

\vspace{0.5cm}

\centerline{\large Kettering University, Flint, MI, USA}

\vspace{0.5cm}

\noindent
{\bf Abstract.} The paper introduces a method of partial fractions with matrix coefficients and its applications to finding chains of generalized eigenvectors, to evaluation of matrix exponentials, and to solution of  linear systems of ordinary differential equations with constant coefficients. 

\noindent
{\bf Keywords:} partial fractions with matrix coefficients, chains of  generalized eigenvectore, systems of differential equations.

\noindent
{\bf MSC:} 34-01, 15-01

\section{Introduction}

In textbooks of ordinary differential equations the Laplace Transform and partial fractions are used for solution of initial value problems for ordinary linear differential equations with constant coefficients.  At the same time, for evaluation of matrix exponentials and for solution of initial value problems for linear systems of ordinary differential equations with constant coefficients  eigenvalues and chains of generalized eigenvectors are used. 

In this paper a different approach is introduced. It is based on representation of the resolvent of a square matrix by the sum of partial fractions with matrix coefficients (see \cite{kato}, p. 40). Combined with the Laplace Transform it presents a more simple technique for solving problems for linear systems of ordinary differential equations with constant coefficients. This procedure is similar to partial fractions students learn in Calculus with the only difference that the coefficients are now matrices. Having matrices as undefined coefficients does not even change the procedure of their evaluation and does not make it more difficult. As it is shown in the paper, the columns of the matrix coefficients of partial fractions are eigenvectors and generalized eigenvectors. Thus, the method of partial fractions with matrix coefficients gives a simple procedure for finding chains of generalized eigenvectors. 

\section{Partial fractions with matrix coefficients and chains of generalized eigenvectors}

In this section the relationship between matrix coefficients of partial fractions and generalized eigenvectors is described. Let $A$ be a square matrix with complex entries and eigenvalues $\l_1$, $\l_2$, \dots, $\l_k$ of the algebraic multiplicities $r_1$, $r_2$, \dots, $r_k$ respectively. Thus,
$$
\det(sI-A)=(s-\l_1)^{r_1}(s-\l_2)^{r_2}\dots(s-\l_k)^{r_k},
$$
$r_1+r_2+\dots+r_k=n$. Then,
$$
(sI-A)^{-1}=\frac{1}{(s-\l_1)^{r_1}(s-\l_2)^{r_2}\dots(s-\l_k)^{r_k}}M(s),
$$
where $M(s)$ is the transposed matrix of cofactors of the matrix $sI-A$ and $I$ is the identity matrix.
There exists a partial fraction decomposition with matrix coefficients of the resolvent $(sI-A)^{-1}$:
\bel{PFDWMC}
(sI-A)^{-1}=\sum_{i=1}^k\sum_{j=1}^{r_i}\frac{1}{(s-\l_i)^j}B_{ij},
\ee
where $B_{ij}$ are matrices that can be determined from this identity  (\cite{kato}, p.40). 
Afer applying $sI-A$ to both parts of this identity one gets
\bel{myeig1}
\sum_{i=1}^k\sum_{j=1}^{r_i}\frac{1}{(s-\l_i)^j}(sI-A)B_{ij}=I.
\ee
Then, multiplication by $(s-\l_1)^{r_1}$ yields
 \bel{myeig2}
\sum_{j=1}^{r_1}(s-\l_1)^{r_1-j}(sI-A)B_{1j}+\sum_{i=2}^k\sum_{j=1}^{r_i}\frac{(s-\l_1)^{r_1}}{(s-\l_i)^j}(sI-A)B_{ij}=(s-\l_1)^{r_1}I.
\ee
After choosing $s=\l_1$, one gets
\bel{myeig0}
(\l_1I-A)B_{1r_1}=0
\ee
or
$$
AB_{1r_1}=\l_1B_{1r_1}.
$$
Let vectors $\bb_1$, $\bb_2$, \dots , $\bb_n$ be the columns of the matrix $B_{1r_1}$, then
$$
A\bb_j=\l_1\bb_j, \hbox{ for } j=1,2,\dots,n.
$$ 
Thus, every nonzero column of the matrix $B_{1r_1}$ is an eigenvector  corresponding to the eigenvalue $\l_1$.
If $r_1>1$, then $m$ differentiations of (\ref{myeig2}) with respect to $s$, $m< r_1$, imply
$$
\sum_{j=1}^{r_1-m}\f{(r_1-j)!}{(r_1-j-m)!}(s-\l_1)^{r_1-j-m}(sI-A)B_{1j}
$$
$$
+m\sum_{j=1}^{r_1-m+1}\f{(r_1-j)!}{(r_1-j-m+1)!}(s-\l_1)^{r_1-j-m+1}B_{1j}
$$
$$
=\f{d^m}{ds^m}\left((s-\l_1)^{r_1}\left(I-\sum_{i=2}^k\sum_{j=1}^{r_i}\frac{1}{(s-\l_i)^j}(sI-A)B_{ij}\right)\right).
$$
Since, $m<r_1$, the substitution $s=\l_1$ yields
$$
m!(\l_1 I-A)B_{1,r_1-m}+m!B_{1,r_1-m+1}=0.
$$
Therefore,
\bel{myeig3}
(A-\l_1I)B_{1,j}=B_{1,j+1},\q j=1,2,\dots,r_1-1,
\ee
and
\bel{myeig31}
B_{1,j}=(A-\l_1I)^{j-1}B_{1,1},\q j=1,2,\dots,r_1.
\ee
Let $\bb_{j1}$, $\bb_{j2}$, \dots, $\bb_{jn}$ be the columns of the matrix $B_{1j}$, that is, $B_{1j}=[\bb_{j1}, \bb_{j2}, \dots, \bb_{jn}]$. For $1\le m\le n$ consider the sequence $\bb_{1m}$, $\bb_{2m}$, \dots, $\bb_{r_1m}$ of the $m$th columns of the matrices $B_{11}$, $B_{12}$, \dots, $B_{1r_1}$. Formulas (\ref{myeig31}) imply
\bel{myeig32}
(A-\l_1I)\bb_{j,m}=\bb_{j+1,m},\q j=1,2,\dots,r_1-1.
\ee
Thus, if $\bb_{jm}=\0$ for some $j$, then $\bb_{im}=\0$ for $j<i\le r_i$.
Let $l_m=\max\{j\,|\, \bb_{jm}\neq\0\}$. Then, nonzero vectors $\bb_{1m}$, $\bb_{2m}$, \dots, $\bb_{l_mm}$ form a chain of vectors such that
\bel{myeig33}
(A-\l_1I)\bb_{l_m,m}=\0, (A-\l_1I)\bb_{l_m-1,m}=\bb_{l_m,m},\dots,  (A-\l_1I)\bb_{1m}=\bb_{2m}.
\ee
Therefore, $\bb_{l_m,m}$ is an eigenvector and $\bb_{j,m}$ are generalized eigenvectors of ranks $l_m-j+1$ respectively, $j=1,\dots,l_{m}-1$. Thus, the nonzero columns $\bb_{1m}$, $\bb_{2m}$, \dots, $\bb_{l_mm}$ form a chain of generalized eigenvectors $\bb_{1m}$, $\bb_{2m}$, \dots, $\bb_{l_mm}$ of the length $l_m$ of the ranks $l_m$, $l_m-1$,\dots,$1$. 
Applying similar procedure to the remaining eigenvalues one gets the following theorem.

\begin{theorem} 
Let $A$ be an $n\times n$ matrix with complex entries and eigenvalues $\l_1$, $\l_2$, \dots, $\l_k$ of the algebraic multiplicities $r_1$, $r_2$, \dots, $r_k$ respectively. 

(i) There exists a unique partial fraction decomposition of the resolvent of $A$
\bel{t1}
(sI-A)^{-1}=\sum_{i=1}^k\sum_{j=1}^{r_i}\frac{1}{(s-\l_i)^j}B_{ij}.
\ee

(ii) If $\bb$  is a nonzero column of the matrix $B_{ir_i}$, then $\bb$ is an eigenvector of $A$ corresponding to the eigenvalue $\l_i$.

(iii) If $\bb$  is a nonzero column of the matrix $B_{ij}$, then $\bb$ is a generalized eigenvector of $A$ of the rank less than or equal to $r_i+1-j$ corresponding to the eigenvalue $\l_i$.

(iv) For $i=1,\dots,k$ and $m=1,\dots,n$ nonzero $m$th columns of the matrices $B_{ij}$, $j=1,\dots,r_i$,  form a chain of generalized eigenvectors corresponding to $\l_i$ ending with an eigenvector.

(v) If $\v$ is a generalized eigenvector of $A$ corresponding to $\l_i$, then $\v$ belongs to the column space of each of the matrices $B_{i1}$, $B_{i2}$, \dots, $B_{ij_0}$, where \break $j_0=1+\max\{m\, |\, \v\in Image((A-\l_iI)^m)\}$.
\end{theorem}

\noindent{\bf Proof. }

\noindent (i) The partial fraction decomposition (\ref{PFDWMC}) can be considered as $n^2$ partial fraction decompositions applied to the entries of the matrix $(sI-A)^{-1}$. Thus, existence and uniqueness of such decomposition follows immediately from the existence and uniqueness in the scalar case. 

\noindent (ii) follows from (\ref{myeig0}).

\noindent (iii) follows from (\ref{myeig0}) and (\ref{myeig3}), since they imply that $(A-\l_i I)^{r_i+1-j}B_{ij}=\0$.

\noindent (iv) follows from  (\ref{myeig33}).

\noindent (v) Let $V_p$ be the invariant subspace corresponding to the eigenvalue $\l_p$ and $\Gamma$ be a positively oriented contour in the complex plane enclosing $s=\l_p$ and excluding all other eigenvalues of $A$ . Then, (see \cite{kato}),
\bel{fkato}
\f{1}{2\pi i}\int\limits_\Gamma(sI-A)^{-1}ds=\sum_{i=1}^k\sum_{j=1}^{r_i}\left(\f{1}{2\pi i}\int\limits_\Gamma\frac{1}{(s-\l_i)^j}ds\right)B_{ij}=B_{p1}
\ee
and, therefore,  $B_{p1}$ is a projector on $V_p$. Denote by $V_{pj}$ the column space of the matrix $B_{pj}$. Then $V_{p1}=V_p$ contains all generalized eigenvectors corresponding to $\l_p$. Hence, if $\v$ is a generalized eigenvector 
corresponding to $\l_i$ then $\v\in V_{i1}$.
The assumption $\v\in Image((A-\l_iI)^{j_0-1})$ implies that there exist generalized eigenvectors $\v_1$, $\v_2$, \dots, $\v_{j_0-1}$ in $V_{i1}$ such that $\v=(A-\l_iI)^{j}\v_j$. Formula (\ref{myeig31}) implies that $(A-\l_iI)^j$ maps the column space of the matrix $B_{i1}$ into the column space of $B_{i,j+1}$. Thus, $\v$ belongs to the column spaces of the matrices $B_{i1}$, $B_{i2}$, \dots, $B_{ij_0}$.
\hspace{10.3cm}$\square$

\begin{remark}
Formula (\ref{fkato}) implies that $(A-\l_i I)B_{i1}=B_{i1}(A-\l_i I)$ and $B_{i1}^2=B_{i1}$. Then, $B_{ij}=(A-\l_i)^{j-1}B_{i1}=(A-\l_i)^{j-1}B_{i1}^{j-1}=((A-\l_i)B_{i1})^{j-1}=B_{i2}^{j-1}$ for $j=2,\dots,r_i$.
\end{remark}
\begin{example}
Find the chains of generalized eigenvectors of the matrix of \break $A=\begin{bmatrix}0&1&2\cr -2&4&0\cr -1&1&2\end{bmatrix}$.
\end{example}
First,
$$
\det(sI-A)=\begin{vmatrix}s&-1&-2\cr 2&s-4&0\cr 1&-1&s-2\end{vmatrix}=s^3-6s^2+12s-8=(s-2)^3.
$$
Thus, the matrix $A$ has an eigenvalue $\l=2$ of the algebraic multiplicity 3. Therefore, the partial fraction decomposition of the resolvent is
$$
(sI-A)^{-1}=\frac{1}{s-2}B_1+\frac{1}{(s-2)^2}B_2+\frac{1}{(s-2)^3}B_3.
$$
Then,
$$
(sI-A)^{-1}=\frac{1}{(s-2)^3}M(s),
$$
where $M(s)$ is the transposed matrix of the cofactors of $sI-A$,
$$
M(s)=\begin{bmatrix}s^2-6s+8&s&2s-8\cr -2s+4&s^2-2s+2&-4\cr -s+2&s-1&s^2-4s+2\end{bmatrix}.
$$
Multiplication by $(s-2)^3$ yields the identity
$$
(s-2)^2B_1+(s-2)B_2+B_3=\begin{bmatrix}s^2-6s+8&s&2s-8\cr -2s+4&s^2-2s+2&-4\cr -s+2&s-1&s^2-4s+2\end{bmatrix}.
$$ 

Thus,
$$
s^2B_1+s(B_2-4B_1)+B_3-2B_2+4B_1
$$
$$
=s^2\begin{bmatrix}1&0&0\cr 0&1&0\cr 0&0&1\end{bmatrix}
+s\begin{bmatrix}-6&1&2\cr -2&-2&0\cr -1&1&-4\end{bmatrix}+\begin{bmatrix}8&0&-8\cr 4&2&-4\cr 2&-1&2\end{bmatrix}.
$$
Therefore,
$$
B_1=\begin{bmatrix}1&0&0\cr 0&1&0\cr 0&0&1\end{bmatrix},\q B_2=4B_1+\begin{bmatrix}-6&1&2\cr -2&-2&0\cr -1&1&-4\end{bmatrix}=\begin{bmatrix}-2&1&2\cr -2&2&0\cr -1&1&0\end{bmatrix},
$$
$$
B_3=2B_2-4B_1+\begin{bmatrix}8&0&-8\cr 4&2&-4\cr 2&-1&2\end{bmatrix}=\begin{bmatrix}0&2&-4\cr 0&2&-4\cr 0&1&-2\end{bmatrix}.
$$
where
$$
B_1=\begin{bmatrix}1&0&0\cr 0&1&0\cr 0&0&1\end{bmatrix},\q B_2=\begin{bmatrix}-2&1&2\cr -2&2&0\cr -1&1&0\end{bmatrix},
B_3=\begin{bmatrix}0&2&-4\cr 0&2&-4\cr 0&1&-2\end{bmatrix}.
$$
Then, as it was shown above in (\ref{myeig0}), (\ref{myeig3}),
$$
(A-2I)B_3=0\q (A-2I)B_2=B_3\q (A-2I)B_1=B_2
$$
Therefore, after choosing the first, the second, and the third columns one gets the following chains of generalized eigenvectors ending with an eigenvector:
$$
\v_1=\begin{bmatrix}1\cr 0\cr 0\end{bmatrix},\q \v_2=\begin{bmatrix}-2\cr -2\cr -1\end{bmatrix},
$$
or
$$
\v_1=\begin{bmatrix}0\cr 1\cr 0\end{bmatrix},\q \v_2=\begin{bmatrix}1\cr 2\cr 1\end{bmatrix},\q \v_3=\begin{bmatrix}2\cr 2\cr 1\end{bmatrix},
$$
or
$$
\v_1=\begin{bmatrix}0\cr 0\cr 1\end{bmatrix},\q \v_2=\begin{bmatrix}2\cr 0\cr 0\end{bmatrix},\q \v_3=\begin{bmatrix}-4\cr -4\cr -2\end{bmatrix}.
$$

\section{Matrix exponential and  linear systems with constant coefficients}

In this section the partial fraction decomposition is applied to evaluation of the matrix exponential and to solution of a system of linear ordinary differential equations with constant coefficients. 

Let $A$ be an $n\times n$ matrix.  Solution of the initial value problem for the system of differential equations $\dot{\y}=A\y+\fb$ requires evaluation of the matrix exponential $e^{tA}$. Such evaluation may be based on solution of the corresponding eigenvalue problem or on the Laplace Transform. In this paper evaluation of the matrix exponential  $e^{tA}=\L^{-1}\{(sI-A)^{-1}\}$ using the Inverse Laplace transform of the resolvent is discussed. 

For evaluation of the matrix exponential $e^{tA}$ of a matrix $A$ with real entries the characteristic polynomial $\det(sI-A)$ is factored over the field of real numbers. To present the method a typical example is considered, let
$$
\det(sI-A)=(s-a_1)(s-a_2)^3(s^2+a_3^2)[(s+a_4)^2+a_5].
$$
The partial fraction decomposition with undetermined matrix coefficients yields:
$$
(sI-A)^{-1}=\frac{1}{(s-a_1)(s-a_2)^3(s^2+a_3^2)[(s+a_4)^2+\b_4^2]}M(s)
$$
$$
=\frac{1}{s-a_1}B_1+\frac{1}{s-a_2}B_2+\frac{1}{(s-a_2)^2}B_3+\frac{1}{(s-a_2)^3}B_4
$$
$$
+\frac{1}{s^2+a_3^2}(sB_5+B_6)+\frac{1}{(s+a_4)^2+\b_4^2}[(s+a_4)B_7+\b_4 B_8],
$$
where $M(s)$ is the transposed matrix of cofactors of $sI-A$ and $B_1$, $B_2$, \dots , $B_8$ are matrices that can be determined from the identity above as it is shown in examples  below.  Hence,  evaluation of the matrix exponential consists the following steps:\\

(i) factoring $\det(sI-A)$ over the field of real numbers,

(ii) the partial fractions decomposition with undetermined matrix coefficients,

(iii) evaluation of the matrix coefficiens,

(iv) the Inverse Laplace transform.\\

The examples below illustrate this procedure.

\begin{example}
Find the matrix exponential $e^{tA}$ of $A=\begin{bmatrix}6&4\cr -3&-1\end{bmatrix}$.
\end{example}

\noindent Since $sI-A=\begin{bmatrix}s-6&-4\cr 3&s+1\end{bmatrix}$ and $\det(sI-A)=(s-2)(s-3)$,
$$
(sI-A)^{-1}=\frac{1}{(s-2)(s-3)}\begin{bmatrix}s+1&4\cr&\cr-3&s-6\end{bmatrix}
=\frac{1}{s-2}B_1+\frac{1}{s-3}B_2,
$$
where $B_1$ and $B_2$ are $2\times 2$ matrices that have to be found. Multiplication by $(s-2)(s-3)$ yields the identity
$$
(s-3)B_1+(s-2)B_2=\begin{bmatrix}s+1&4\cr&\cr-3&s-6\end{bmatrix}.
$$ 
Substitutions $s=2$ and $s=3$ yield 
$$
-B_1=\begin{bmatrix}3&4\cr&\cr-3&-4\end{bmatrix}\q\hbox{and}\q B_2=\begin{bmatrix}4&4\cr&\cr-3&-3\end{bmatrix}.
$$
Therefore,
$$
e^{tA}=\L^{-1}\{(sI-A)^{-1}\}=\L^{-1}\left\{\frac{1}{s-2}B_1+\frac{1}{s-3}B_2\right\}
$$
$$
=\L^{-1}\left\{\frac{1}{s-2}\right\}B_1+\L^{-1}\left\{\frac{1}{s-3}\right\}B_2=e^{2t}B_1+e^{3t}B_2
$$
$$
=e^{2t}\begin{bmatrix}-3&-4\cr&\cr3&4\end{bmatrix}+e^{3t}\begin{bmatrix}4&4\cr&\cr-3&-3\end{bmatrix}.
$$

\begin{example}
Find the matrix exponential $e^{tA}$ of $A=\begin{bmatrix}5&17\cr -2&-5\end{bmatrix}$.
\end{example}

\noindent In this example $sI-A=\begin{bmatrix}s-5&-17\cr 2&s+5\end{bmatrix}$, $\det(sI-A)=s^2+9$, and
$$
(sI-A)^{-1}=\frac{1}{s^2+9}\begin{bmatrix}s+5&17\cr&\cr-2&s-5\end{bmatrix}
=\frac{s}{s^2+9}B_1+\frac{1}{s^2+9}B_2.
$$
Multiplication by $s^2+9$ yields the identity
$$
sB_1+B_2=\begin{bmatrix}s+5&17\cr&\cr-2&s-5\end{bmatrix}=s\begin{bmatrix}1&0\cr&\cr 0&1\end{bmatrix}+\begin{bmatrix}5&17\cr&\cr-2&-5\end{bmatrix}.
$$ 
Thus,
$$
B_1=\begin{bmatrix}1&0\cr&\cr 0&1\end{bmatrix}\q\hbox{and}\q B_2=\begin{bmatrix}5&17\cr&\cr-2&-5\end{bmatrix}.
$$
Therefore,
$$
e^{tA}=\L^{-1}\{(sI-A)^{-1}\}=\L^{-1}\left\{\frac{s}{s^2+9}B_1+\frac{1}{s^2+9}B_2\right\}
$$
$$
=\L^{-1}\left\{\frac{s}{s^2+9}\right\}B_1+\L^{-1}\left\{\frac{1}{s^2+9}\right\}B_2=\cos(3t)B_1+\frac{1}{3}\sin(3t)B_2
$$
$$
=\cos(3t)\begin{bmatrix}1&0\cr&\cr 0&1\end{bmatrix}+\sin(3t)\begin{bmatrix}\frac{5}{3}&\frac{17}{3}\cr&\cr-\frac{2}{3}&-\frac{5}{3}\end{bmatrix}.
$$

\begin{example}Solve the initial value problem
$$
\y'=\begin{bmatrix}-5&6&2\cr -6&7&2\cr 6&-6&-1\end{bmatrix}\y,\q \y(0)=\begin{bmatrix}1\cr-1\cr 2\end{bmatrix}.
$$
 \end{example}
The characteristic polynomial is
$$
L(s)=(s-1)(s^2-1)=(s-1)^2(s+1).
$$
and
$$
(sI-A)^{-1}=\f{1}{(s-1)^2(s+1)}\begin{bmatrix}s^2-6s+5&6s-6&2s-2\cr -(6s-6)&s^2+6s-7&2s-2\cr 6s-6&-(6s-6)&s^2-2s+1\end{bmatrix}
$$
$$
=\f{1}{s-1}B_1+\f{1}{(s-1)^2}B_2+\f{1}{s+1}B_3.
$$
Therefore
$$
\begin{bmatrix}s^2-6s+5&6s-6&2s-2\cr -(6s-6)&s^2+6s-7&2s-2\cr 6s-6&-(6s-6)&s^2-2s+1\end{bmatrix}=(s-1)(s+1)B_1+(s+1)B_2+(s-1)^2B_3.
$$
Choosing $s=1$, $s=-1$, and $s=0$ yields
$$
B_2=\begin{bmatrix}0&0&0\cr 0&0&0\cr 0&0&0\end{bmatrix},\q B_3=\begin{bmatrix}3&-3&-1\cr 3&-3&-1\cr -3&3&1\end{bmatrix},
$$
$$
-B_1+B_2+B_3=\begin{bmatrix}5&-6&-2\cr 6&-7&-2\cr -6&6&1\end{bmatrix}.
$$
Hence,
$$
B_1=B_2+B_3-\begin{bmatrix}5&-6&-2\cr 6&-7&-2\cr -6&6&1\end{bmatrix}=\begin{bmatrix}-2&3&1\cr -3&4&1\cr 3&-3&0\end{bmatrix}.
$$
Thus,
$$
e^{tA}=\L^{-1}\left\{\f{1}{s-1}\right\}B_1+\L^{-1}\left\{\f{1}{(s-1)^2}\right\}B_2+\L^{-1}\left\{\f{1}{s+1}\right\}B_3
$$
$$
=e^{t}\begin{bmatrix}-2&3&1\cr -3&4&1\cr 3&-3&0\end{bmatrix}+e^{-t}\begin{bmatrix}3&-3&-1\cr 3&-3&-1\cr -3&3&1\end{bmatrix}
$$
and
$$
\y=e^{tA}\y(0)=e^{t}\begin{bmatrix}-2&3&1\cr -3&4&1\cr 3&-3&0\end{bmatrix}\begin{bmatrix}1\cr-1\cr 2\end{bmatrix}+e^{-t}\begin{bmatrix}3&-3&-1\cr 3&-3&-1\cr -3&3&1\end{bmatrix}\begin{bmatrix}1\cr-1\cr 2\end{bmatrix}
$$
$$
=e^{t}\begin{bmatrix}-3\cr-5\cr 6\end{bmatrix}+e^{-t}\begin{bmatrix}4\cr 4 \cr -4\end{bmatrix}.
$$

\begin{example} Find the general solution of
$$
\y'=\begin{bmatrix}1&9&6\cr -6&-20&-12\cr 9&24&13\end{bmatrix}\y.
$$
\end{example}
Thus,
$$
 sI-A=\begin{bmatrix}s-1&-9&-6\cr 6&s+20&12\cr -9&-24&s-13\end{bmatrix}
$$
and the characteristic polynomial is
$$
L(s)=\det(sI-A)=s^3+6s^2+21s+26=(s+2)(s^2+4s+13).
$$
Therefore,
$$
(sI-A)^{-1}=\f{1}{(s+2)[(s+2)^2+9]}\begin{bmatrix}s^2+7s+28&9s+27&6s+12\cr -6s-30&s^2-14s-41&-12s-24\cr 9s+36&24s+57&s^2+19s+34\end{bmatrix},
$$
$$
=\f{1}{s+2}B_1+\f{s+2}{(s+2)^2+9}B_2+\f{3}{(s+2)^2+9}B_3.
$$
This implies the identity
$$
\begin{bmatrix}s^2+7s+28&9s+27&6s+12\cr -6s-30&s^2-14s-41&-12s-24\cr 9s+36&24s+57&s^2+19s+34\end{bmatrix}
$$
$$
=[(s+2)^2+9]B_1+(s+2)^2B_2+3(s+2)B_3.
$$
Choosing $s=-2$, $s=0$, and setting coeficients of $s^2$ equal yields
$$
9B_1=\begin{bmatrix}18&9&0\cr -18&-9&0\cr 18&9&0\end{bmatrix},\q 13B_1+4B_2+6B_3=\begin{bmatrix}28&27&12\cr -30&-41&-24\cr 36&57&34\end{bmatrix},
$$
$$
B_1+B_2=\begin{bmatrix}1&0&0\cr 0&1&0\cr 0&0&1\end{bmatrix}.
$$
Thus,
$$
B_1=\begin{bmatrix}2&1&0\cr -2&-1&0\cr 2&1&0\end{bmatrix},\q B_2=\begin{bmatrix}-1&-1&0\cr 2&2&0\cr -2&-1&1\end{bmatrix},\q B_3=\begin{bmatrix}1&3&2\cr -2&-6&-4\cr 3&8&5\end{bmatrix}.
$$
Then,
$$
e^{tA}=e^{-2t}\begin{bmatrix}2&1&0\cr -2&-1&0\cr 2&1&0\end{bmatrix}+e^{-2t}\cos(3t)\begin{bmatrix}-1&-1&0\cr 2&2&0\cr -2&-1&1\end{bmatrix}
+e^{-2t}\sin(3t)\begin{bmatrix}1&3&2\cr -2&-6&-4\cr 3&8&5\end{bmatrix}.
$$
and
$$
\y=e^{tA}\begin{bmatrix}C_1\cr C_2\cr C_3\end{bmatrix}=C_1\left(e^{-2t}\begin{bmatrix}2\cr -2\cr 2\end{bmatrix}+e^{-2t}\cos(3t)\begin{bmatrix}-1\cr 2\cr -2\end{bmatrix}
+e^{-2t}\sin(3t)\begin{bmatrix}1\cr -2\cr 3\end{bmatrix}\right)
$$
$$
+C_2\left(e^{-2t}\begin{bmatrix}1\cr -1\cr 1\end{bmatrix}+e^{-2t}\cos(3t)\begin{bmatrix}-1\cr 2\cr -1\end{bmatrix}
+e^{-2t}\sin(3t)\begin{bmatrix}3\cr -6\cr 8\end{bmatrix}\right)
$$
$$
+C_3\left(e^{-2t}\cos(3t)\begin{bmatrix}0\cr 0\cr -1\end{bmatrix}
+e^{-2t}\sin(3t)\begin{bmatrix}2\cr -4\cr 5\end{bmatrix}\right).
$$

{\bf Summary.} The paper introduces an alternative method for solving initial value problems for linear systems of ordinary differential equations (ODEs) with constant coefficients. Traditionally, scalar ODEs use the Laplace Transform and partial fractions, while systems of ODEs rely on matrix exponentials, eigenvalues, and generalized eigenvectors.

The proposed approach combines the Laplace Transform with partial fraction decomposition using matrix coefficients. The method mirrors the standard partial fraction technique from Calculus, but the unknowns are matrices rather than scalars. Despite this change, the evaluation process remains essentially the same.

A key insight is that the columns of these matrix coefficients correspond to eigenvectors and generalized eigenvectors of the system matrix. As a result, this method provides a simpler way to construct the chains of generalized eigenvectors, as well as offering a more accessible technique for solving linear systems of ODEs. Thus, it gives an easier and shorter path to applications such as the Control Theory for undergraduate students, especially for students majoring in engineering.\\

\noindent {\bf Conflict of Interest:} The author declares that there is no conflict of interest.

\bibliography{MatrixPFD2}

\begin{thebibliography}{10}

\bibitem{kato}
Kato T.
\newblock Perturbation Theory for Linear Operators.
\newblock Springer; 1995.
\end{thebibliography}

\end{document}